\RequirePackage{etoolbox}

\documentclass[ba]{imsart}
\pubyear{0000}
\volume{00}
\issue{0}
\doi{0000}
\firstpage{1}
\lastpage{14}

\usepackage{amsthm}
\usepackage{amsmath}
\usepackage{natbib}
\usepackage[colorlinks,citecolor=blue,urlcolor=blue,filecolor=blue,backref=page]{hyperref}
\usepackage{graphicx}

\startlocaldefs
\numberwithin{equation}{section}
\theoremstyle{plain}

\endlocaldefs
\usepackage{amsmath}
\usepackage{graphicx}
\usepackage{natbib}
\usepackage[colorlinks,citecolor=blue,urlcolor=blue,filecolor=blue,backref=page]{hyperref}

\begin{document}
\begin{frontmatter}
\title{Optimal proposals for Approximate Bayesian Computation\thanksref{T1}}
\thankstext{T1}{Work completed at the Center for Computational Astrophysics, Flatiron Institute, 162 5th Ave, New York City, NY 10010, USA}
\runtitle{Optimal proposals for ABC}
\begin{aug}

\author{\fnms{Justin} \snm{Alsing}\thanksref{a1}\ead[label=e1]{jalsing@flatironinstitute.org}}
\author{\fnms{Benjamin D.} \snm{Wandelt}\thanksref{a1}\thanksref{a2}\thanksref{a3}\thanksref{a4}\ead[label=e2]{bwandelt@flatironinstitute.org}}
\author{\fnms{Stephen M.} \snm{Feeney}\thanksref{a1}\ead[label=e3]{sfeeney@flatironinstitute.org}}
  
\address[a1]{Center for Computational Astrophysics, Flatiron Institute, 162 5th Ave, New York City, NY 10010, USA}
\address[a2]{Sorbonne Universit\'es, CNRS, UMR 7095, Institut d'Astrophysique de Paris, 98 \textit{bis} boulevard Arago, 75014 Paris, France}
\address[a3]{Sorbonne Universit\'es, Institut  Lagrange  de  Paris  (ILP),  98 \textit{bis} boulevard Arago, 75014 Paris, France}
\address[a4]{Department of Astrophysical Sciences, Princeton University, Princeton, NJ 08540, USA}
\runauthor{ J. Alsing, B. Wandelt \& S. Feeney}
\end{aug}

\begin{abstract}
We derive the optimal proposal density for Approximate Bayesian Computation (ABC) using Sequential Monte Carlo (SMC) (or Population Monte Carlo, PMC). The criterion for optimality is that the SMC/PMC-ABC sampler maximise the effective number of samples per parameter proposal. The optimal proposal density represents the optimal trade-off  between favoring high acceptance rate and reducing the variance of the importance weights of accepted samples. We discuss two convenient approximations of this proposal and show that the optimal proposal density gives a significant boost in the expected sampling efficiency compared to standard kernels that are in common use in the ABC literature, especially as the number of parameters increases.
\end{abstract}


\end{frontmatter}

\section{Introduction}
The na\"ive approach to Approximate Bayesian Computation (ABC) generates (compressed) data simulations for parameters $\theta$ that are drawn from the prior $\pi(\theta)$. If the resulting simulated data $d'$ is within $\epsilon$ of the true data $d$, i.e. under a distance metric $\rho(d,d')\leq\epsilon$, then $\theta$ is accepted as a sample from the approximate posterior density, $p(\theta | \rho(d,d')\leq\epsilon)$.  

If simulation is costly, it is advantageous to attempt to increase the fraction of accepted $\theta$ by proposing new candidate $\theta$s from a proposal density $q(\theta)$ that is large in parameter ranges that are preferred by the data, and small in less interesting regions of the prior volume. It is natural to base the choice of $q(\theta)$ on the current accepted samples from the approximate posterior $p(\theta | \rho(d,d')\leq\epsilon)$. This is the approach taken in SMC-ABC algorithms (see eg., \citealp{ABCreview2012} for a review); a proposal density $q(\theta)$ is constructed at each population iteration, which is typically a kernel density estimator (KDE) based on the accepted points from the previous population (ie, a kernel that adapts as the algorithm steps through successive population iterations \citealp{Beaumont2009, Mckinley2009, Toni2009a, Barnes2011, Didelot2011, Jasra2012, Filippi2013, Bonassi2015}).

The price to pay for the increased fraction of accepted ABC samples of $\theta$ is the necessity to importance weight the accepted samples by $\pi/q$. The variance in these importance weights will reduce the effective number of samples. The more concentrated $q$ is relative to $\pi$, on parameters $\theta$ that have a high probability of being accepted, the larger the variance in the importance weights $\pi/q$. In practice, a poor choice of proposal density $q$ can lead to a proposed sample whose proposal density $q$ was low, but is subsequently accepted, leading to a large importance weight that can overwhelm the rest of the (weighted) samples leading to a very small effective sample size. The natural question, then, is how to choose the proposal density that represents the optimal trade-off between a high acceptance rate and a low importance-weight variance.

Most SMC-ABC implementations propose new parameters for forward simulation via an importance weighted KDE based on the previous population's accepted samples, with uniform \citep{Toni2009,Toni2009a,Mckinley2009}, student-t \citep{Didelot2011} and Gaussian kernels \citep{Sisson2007,Beaumont2009,Filippi2013,Bonassi2015} in common use. A choice has to be made for the kernel bandwidth and there are various choices in the literature: for example, for Gaussian kernels, \cite{Sisson2007} use the importance-weighted variance of the previous population samples, \citet{Beaumont2009} use twice the importance-weighted variance of the previous population samples, whilst \cite{Bonassi2015} use standard recommendations from \citet{West1993} and \citet{Scott2005}\footnote{ie., taking the (component-wise) importance-weighted variance divided by $N^{1/3}$, for $N$ samples.}. Previous studies have sought kernels that are optimal in the following sense: they minimize the sum of the Kullback-Leibler divergence between the proposal KDE and the target density, and the negative log acceptance ratio, providing some trade-off between closeness of the target and proposal and the acceptance ratio \citep{Beaumont2009,Filippi2013}. \citet{Beaumont2009} showed that for a global Gaussian kernel, this optimality criterion leads to a bandwidth equal to twice the importance-weighted variance of the previous population's accepted samples (with some further refinement and generalization by \citealp{Filippi2013}). \citet{Filippi2013} also considered the same optimality criterion applied to local rather than global kernels, deriving a locally-optimal kernel-covariance. Whilst this optimality criterion has proved powerful, the relative importance and utility of the KL divergence and acceptance ratio terms in this approach is ambiguous.

In this paper we take a slightly different approach and derive the proposal density $q(\theta)$ for SMC-ABC sampling that maximizes the effective number of samples per parameter proposal (and hence forward simulation). This provides the optimal trade-off between high acceptance rate and low variance in the importance weights under a straightforward and pragmatic definition of optimality. Rather than restricting to a given class of perturbation kernels (eg., Gaussian kernels), we derive the optimal proposal density (in the asymptotic $\epsilon\rightarrow 0$ limit) assuming only that some density estimator for the ABC posterior is available at each population iteration. The result provides a well-motivated guide for adaptive proposal density choice for SMC-ABC sampling.
\section{Optimal SMC-ABC proposal densities}
\label{sec:optimal}
We define the \textit{sampling efficiency} as the functional $\omega[q]$ that measures the effective number of samples per parameter proposal. This is composed of two components: 1) $f_{a}$, the fraction of proposed points that will be accepted, and 2) $N_w$, the effective number of points after application of the importance weights.

The expected fraction of accepted points is given by
\begin{align}
f_{a}=P(\rho(d,d')<\epsilon)&=\int I_{\rho(d,d')<\epsilon}\, p_q(d') dd'\nonumber \\
&\underset{\epsilon
\rightarrow 0}{\approx} V_\epsilon p_q(d)=V_\epsilon \int p(d|\theta)q(\theta) d\theta =V_\epsilon \, p(d) \int \frac{q(\theta)}{\pi(\theta)}p(\theta|d)d\theta 
\end{align}
where $ p_q(d)$ is the probability density of simulated data $d$ when the parameters are proposed from $q$, $V_\epsilon$ is the volume of the space for accepted samples, and "$\approx$" becomes accurate in the limit of small $\epsilon$, assuming $p(d)$ is continuous.

The expected effective number of points after application of the importance weights to $N_s$ accepted samples is given by
\begin{align}
N_w=\frac{\left(\sum w_i\right)^2}{\sum w_i^2}=
\frac{\left(\sum_i^{N_s} \pi(\theta_i)/q(\theta_i)\right)^2}{\sum_i^{N_s} \left(\pi(\theta_i)/q(\theta_i)\right)^2}&\\
\underset{N_s\rightarrow\infty}{\approx}\frac{N_s\left(\int \frac{\pi}{q} p_\mathrm{accepted}(\theta)d\theta\right)^2}{\int  \frac{\pi^2}{q^2}p_\mathrm{accepted}(\theta)d\theta}&\\
\underset{\epsilon
\rightarrow 0}{\approx}\frac{N_s V_\epsilon\,\left(\int  \frac{\pi}{q} p(d|\theta)q(\theta)d\theta\right)^2}{\int  \frac{\pi^2}{q^2} p(d|\theta)q(\theta)d\theta} &\\
=\frac{N_s V_\epsilon\,p(d)}{\int\frac{\pi(\theta)}{q(\theta)} p(\theta | d) d\theta}. 
\end{align}
We therefore find for the sampling efficiency $\omega[q]$ (dropping $q$-independent constant pre-factors),
\begin{equation}
\label{omega}
\omega[q]\equiv \frac{A[q]}{B[q]} = \frac{\int \frac{q(\theta)}{\pi(\theta)}p(\theta | d)d\theta}{\int \frac{\pi(\theta)}{q(\theta)}p(\theta | d)d\theta}\propto f_a\frac{N_w}{N_s},
\end{equation}
where the second equality defines the functionals $A[q]$ and $B[q]$, and the last proportionality defines the sampling efficiency as the number of effective (accepted) samples per parameter proposal (and hence forward simulation). For na\"{i}ve ABC, proposing parameters from the prior, $\omega[\pi]=1$.
Maximizing $\omega[q]$ with respect to $q$, under the constraint of $q$ being normalized, gives the optimal proposal 
\begin{equation}
\label{eq:optimalq}
q^\ast(\theta)=\sqrt[]{\omega[q^\ast]\frac{p(\theta | d)\pi(\theta)}{2A[q^\ast]-\frac{p(\theta|d)}{\pi(\theta)}}}.
\end{equation}
This is the main result of the paper. This implicit equation for $q^\ast$ can be solved iteratively; alternatively one can just consider $\omega^\ast\equiv\omega[q^\ast]$ and $A^\ast\equiv A[q^\ast]$ as parameters and  
\begin{equation}
q^\ast(\theta)=\max_{\omega^\ast,A^\ast} \omega[q(\theta,\omega^\ast,A^\ast)].
\end{equation}
Further,  there are two simple and fast approximations for $q^\ast$ that can be obtained without iteration or that can inform good starting points for the iteration; these are discussed in \S \ref{sec:approximations}.

Note that we only need to know $A$ to be able to sample from the optimal proposal \eqref{eq:optimalq}, since $\omega$ just sets the normalization. Assuming a density estimator for $p(\theta|d)$ is available, and $\pi$ is easy to evaluate, then Markov Chain Monte Carlo (MCMC) methods can efficiently generate draws from $q^\ast$ to serve as parameter proposals for the next SMC iteration.

The optimal proposal Eq. \eqref{eq:optimalq} can be seen as the geometric mean of posterior and prior (ie., the numerator), with a relative boost where the posterior is larger than the prior. The denominator boosts $q^\ast$ in the region where the posterior peaks, but the peaks in $q^\ast$ are narrower than the corresponding ones in $p(\theta|d)$ in their immediate neighborhoods to compensate for the heavier tails away from the peaks. Illustrative examples are shown in Figures \ref{fig:gaussian_qplots}--\ref{fig:chi2_qplots} and discussed in \S \ref{sec:examples}.

\subsection{Fast approximations to the optimal proposal density}
\label{sec:approximations}
\subsubsection{Geometric mean approximation}
Expanding $q^\ast$ gives the convergent series 
\begin{equation}
q^\ast=\sqrt[]{\frac{p(\theta|d)\pi(\theta)}{2B[q^\ast]}}\sum_{i=0}^{\infty}
\begin{pmatrix}
i-\frac12\\i
\end{pmatrix}\left(\frac{p(\theta|d)}{2A[q^\ast]\pi(\theta)}\right)^i.
\end{equation}
Convergence follows from the properties of the binomial coefficient and the lower bound on $A[q^\ast]$ in Eq.~(\ref{holders}) below. In numerical experiments we find the $i>0$ terms give sub-dominant contributions to $\omega$ for a wide variety of choices for $q$, $p$ and $\pi$. Taking the normalized leading term of the series gives the \textit{geometric mean approximation} to the optimal proposal\footnote{
From Jensen's inequality we have 
\begin{equation}
\int \left(\frac{p}{\pi}\right)^{2}\pi d\theta\geq
\int \left(\frac{p}{\pi}\right)^{\frac32}\pi d\theta\geq
1\geq
\int \left(\frac{p}{\pi}\right)^{\frac12}\pi d\theta.
\end{equation}
It follows that  $\omega[q_0]\geq \omega[\pi]$ and  $\omega[p]\geq \omega[\pi]$  for any pdf $p$ and $\pi$. While we cannot conclude in general that $\omega[q_0]\geq \omega[p]$, Figure \ref{fig:approximate_improvement_posterior} shows that in the Gaussian example $q_0$  outperforms the posterior in all practically relevant cases. },
\begin{equation}
\label{geometric_mean}
q_0 \propto\sqrt{p(\theta | d)\pi(\theta)}.
\end{equation}
Note that for simple cases where the posterior is Gaussian under a uniform prior, the geometric mean approximation leads to a Gaussian proposal density centered on the posterior mean and with twice the posterior covariance. This is similar to taking a Gaussian KDE proposal with bandwidth equal to the estimated posterior variance; this kernel scheme is sometimes adopted in the ABC literature (eg., \citealp{Sisson2007, Ishida2015}).

\subsubsection{Bounded approximation}
Given the value of $A[q^*]$, we could sample the optimal proposal density Eq. \eqref{eq:optimalq}. Whilst $A[q^*]$ is not available a priori, we can bound $A[q^*]$ from above and below: using H\"older's inequality to obtain the upper bound, and the non-negativity of probability densities for the lower bound, we find that 
\begin{equation}
\label{holders}
\frac{1}{2}\sup_\theta\frac{p(\theta|d)}{\pi(\theta)}< A[q^*] \leq\sup_\theta \frac{p(\theta|d)}{\pi(\theta)}.
\end{equation}
Choosing a value of $A[q^*]$ between these bounds avoids the need for iteration for the cost of a mild reduction in optimality. Fixing $A$ to be the average of the upper and the lower bound leads to the \emph{bounded approximation} for the optimal proposal,
\begin{equation}
\label{bounded_approx}
q_{\bar{A}} \propto \sqrt{ \frac{p(\theta | d)\pi(\theta)}{2\bar{A} - \frac{p(\theta | d)}{\pi(\theta)} }},
\end{equation}
with
\begin{equation}
\bar{A} = \frac{3}{4}\sup_\theta\frac{p(\theta|d)}{\pi(\theta)}.
\end{equation}
In numerical experiments, we find that the bounded approximation gives close to optimal sampling efficiencies (see \S \ref{sec:examples} and Table \ref{tab:sampling_efficiencies}). If further optimality is desired, this can be used as a starting point for iterating Eq. \eqref{eq:optimalq} towards the optimal $q^\ast$.

All that is required to propose samples from the optimal proposal or an approximation of it is a density estimator for the posterior. In practice, this could be a KDE or mixture model fit to the accepted samples in the previous SMC population. With a posterior density estimator in hand, the optimal proposal can be found iteratively using Eq. \eqref{eq:optimalq}, or via one of the convenient approximations Eq. \eqref{geometric_mean} or \eqref{bounded_approx}, and then sampled using MCMC or otherwise to generate parameter proposals for the next SMC iteration.


\begin{figure}
\centering
\includegraphics[width=.79\textwidth]{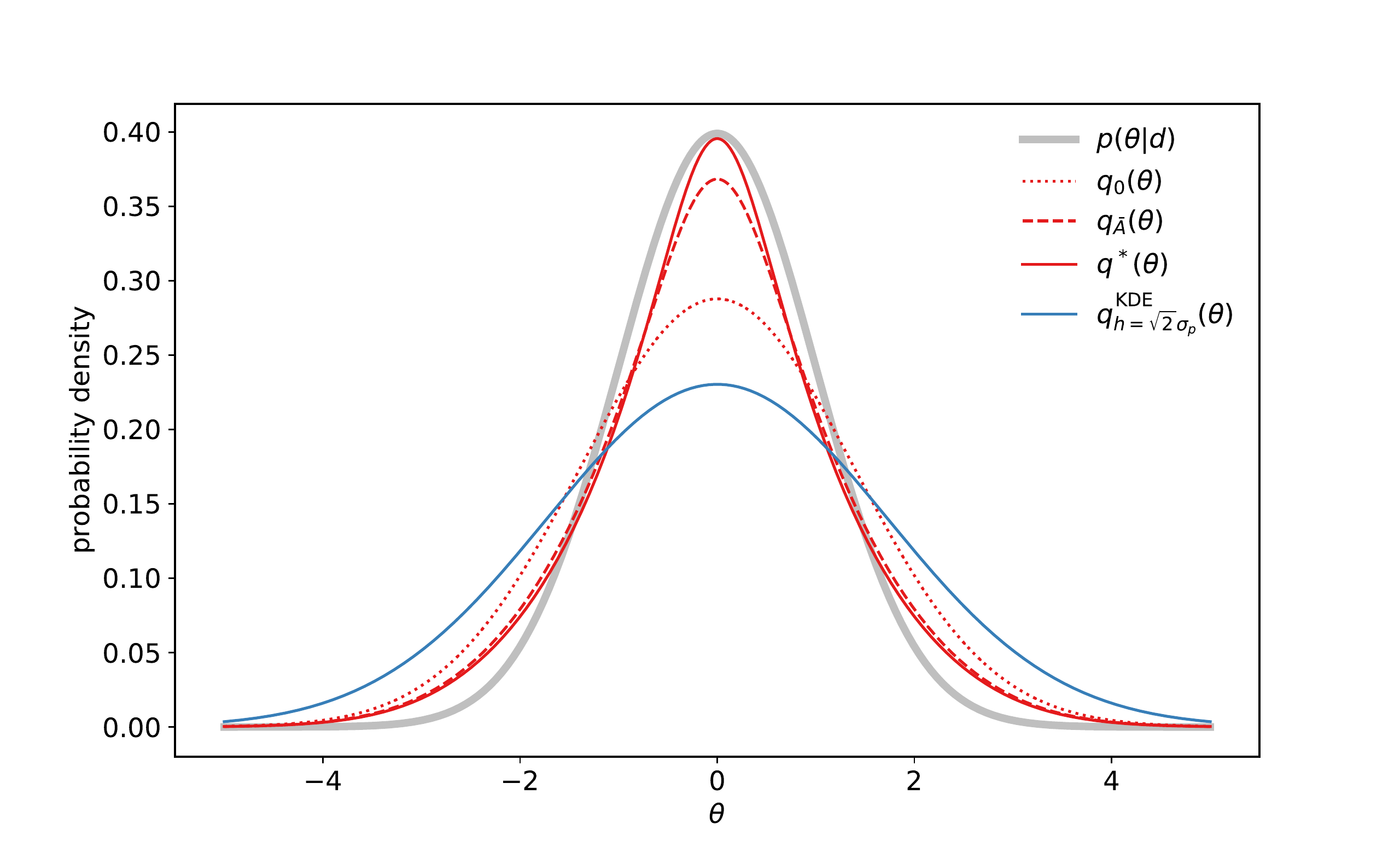}
\caption{ABC proposals for Gaussian posterior $p(\theta | d) = \mathcal{N}(0, 1)$ (grey), with a Gaussian prior $\pi(\theta)=\mathcal{N}(0,5)$ (not shown). From bottom to top at the peak: commonly used KDE proposal with bandwidth of twice the (estimated) posterior variance (blue), geometric mean approximation of the optimal proposal (red-dotted), bounded approximation of the optimal proposal with $A = 3/4\,\sup_\theta p(\theta|d)/\pi(\theta)$ (red-dashed), optimal proposal density (red).}\label{fig:gaussian_qplots}
\end{figure}
\begin{figure}
\centering
\includegraphics[width=.79\textwidth]{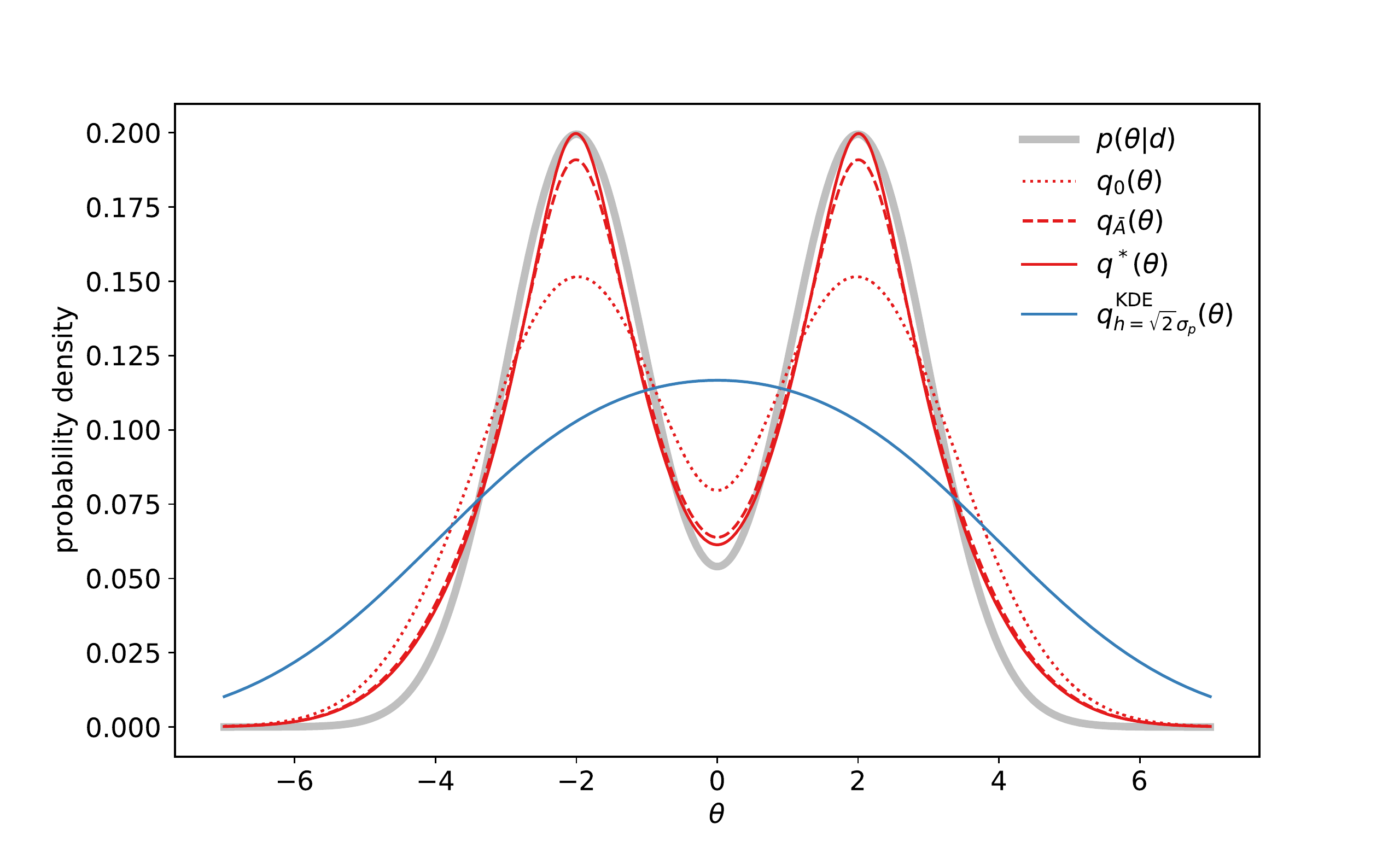}
\caption{ABC proposals for a bimodal posterior $p(\theta | d) = 1/2\,\mathcal{N}(-2, 1) + 1/2\,\mathcal{N}(2, 1)$ (grey), with a Gaussian prior $\pi(\theta)=\mathcal{N}(0,10)$ (not shown). From bottom to top at the peak: commonly used KDE proposal with bandwidth of twice the (estimated) posterior variance (blue), geometric mean approximation of the optimal proposal (red-dotted), bounded approximation of the optimal proposal with $A = 3/4\,\sup_\theta p(\theta|d)/\pi(\theta)$ (red-dashed), optimal proposal density (red).}\label{fig:bimodal_qplots}
\end{figure}
\begin{figure}
\centering
\includegraphics[width=.79\textwidth]{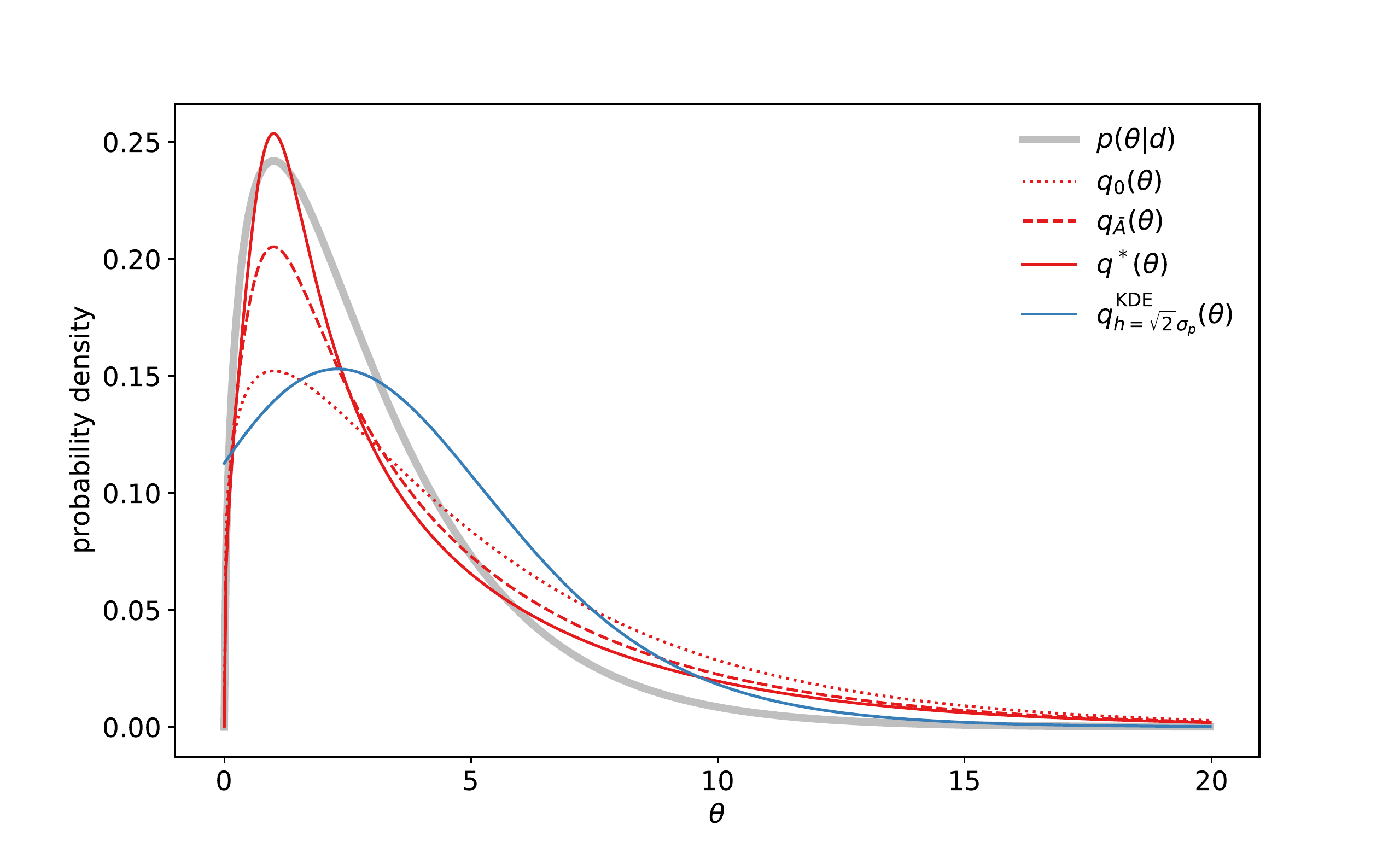}
\caption{ABC proposals for $\chi^2$ posterior $p(\theta | d) = \chi^2(\nu=3)$ (grey), with a uniform prior $\pi(\theta)=\mathcal{U}(0,30)$ (not shown). From bottom to top at the peak: commonly used KDE proposal with bandwidth of twice the (estimated) posterior variance (blue), geometric mean approximation of the optimal proposal (red-dotted), bounded approximation of the optimal proposal with $A = 3/4\,\sup_\theta p(\theta|d)/\pi(\theta)$ (red-dashed), optimal proposal density (red).}\label{fig:chi2_qplots}
\end{figure}
\begin{table*}
\caption{Comparison of the performance of different SMC proposal schemes for the three cases considered. 
Functionals $A[q]$ and $B[q]$ are proportional to the expected acceptance rate and (inverse) importance-weight variance respectively, and $\omega[q]$ defines the sampling efficiency; see Eq. \protect\eqref{omega} for definitions.}
\begin{tabular}{ccccc}
\hline
Proposal, $q(\theta)$ & $A[q]$ & $B[q]$ & $\omega[q]\equiv A[q]/B[q]$ \\
\hline
\multicolumn{1}{l}{{Case I: $p(\theta | d) = \mathcal{N}(0, 1)$, $\pi(\theta)=\mathcal{N}(0,5)$ } }\\
$p(\theta | d)$ &  3.57 & 1.0 & 3.57 \\
$q^\mathrm{KDE}_{h=\sqrt{2}\sigma_p}(\theta)$ & 2.54 & 0.41 & 6.16 \\
$q_0(\theta)$ & 2.96 & 0.38 & 7.71 \\
$q_{\bar{A}}(\theta)$ & 3.26 & 0.40 & 8.20 \\
$q^\ast(\theta)$ & 3.34 & 0.41 & \bf{8.22} \\ \\ 

\multicolumn{2}{l}{Case II: $p(\theta | d) = 1/2\mathcal{N}(-2, 1) + 1/2\mathcal{N}(2, 1)$, $\pi(\theta)=\mathcal{N}(0,10)$}  \\
$p(\theta | d)$ &  3.68 & 1.0 & 3.68 \\
$q^\mathrm{KDE}_{h=\sqrt{2}\sigma_p}(\theta)$  & 2.55 & 0.40 & 6.36 \\
$q_0(\theta)$  & 3.22 & 0.34 & 9.47 \\
$q_{\bar{A}}(\theta)$  & 3.48 & 0.35 & 9.93 \\
$q^\ast(\theta)$  & 3.52 & 0.35 & \bf{9.94} \\ \\

\multicolumn{1}{l}{Case III: $p(\theta | d) = \chi^2(\nu=3)$, $\pi(\theta)=\mathcal{U}(0,30)$}  \\
$p(\theta | d)$  & 4.77 & 1.0 & 4.77 \\
$q^\mathrm{KDE}_{h=\sqrt{2}\sigma_p}(\theta)$  & 3.80 & 0.57 & 6.68 \\
$q_0(\theta)$  & 3.56 & 0.35 & 10.23 \\
$q_{\bar{A}}(\theta)$  & 4.05 & 0.36 & 11.23 \\
$q^\ast(\theta)$  & 4.34 & 0.38 & \bf{11.40} \\

\hline
\end{tabular}
\label{tab:sampling_efficiencies}
\end{table*}

\subsection{Examples}
\label{sec:examples}
\subsubsection{Numerical examples of optimal ABC proposals and their approximations}
Figures \ref{fig:gaussian_qplots}--\ref{fig:chi2_qplots} and Table \ref{tab:sampling_efficiencies} illustrate the optimal ABC proposal density and its approximations in three scenarios: a Gaussian posterior (Figure \ref{fig:gaussian_qplots}), a bimodal double-Gaussian posterior (Figure \ref{fig:bimodal_qplots}), and a $\chi^2$ posterior (Figure \ref{fig:chi2_qplots}). We compare the optimal proposals to the commonly used ABC proposal scheme recommended in \citet{Beaumont2009} for reference: a Gaussian KDE with bandwidth set to double the (estimated) posterior variance. For illustration, the proposals are compared in the converged $\epsilon\rightarrow0$ limit where the approximate posterior (in practice, a density estimator for the accepted samples) is close to the true posterior. The KDE proposal \citep{Beaumont2009} is shown as the convolution of the true posterior with a Gaussian with twice the posterior variance.

The three examples shown in Figures \ref{fig:gaussian_qplots}--\ref{fig:chi2_qplots} and Table \ref{tab:sampling_efficiencies} display the same essential characteristics. The optimal proposals are boosted in regions of high posterior density (around the peak) to give a high acceptance rate, whilst having slightly broader tails compared to the posterior to ensure the importance-weight variance is kept under control (hence giving an improved effective sample size). In contrast, the KDE proposals are typically much broader than the optimal proposals around the posterior peak, which leads to lower expected acceptance rates. Meanwhile, using a posterior density estimate for the proposal gives a poor expected importance-weight variance, owing to the narrower tails compared to the other proposal schemes. The optimal proposal represents the trade-off between high proposal density in regions of high posterior density, and fatter tails to maintain a lower importance-weight variance.

In all three examples shown, the bounded approximation for the optimal proposal Eq. \eqref{bounded_approx} performs nearly as well as the optimal proposal. This is especially clear from Table \ref{tab:sampling_efficiencies}, where the sampling efficiencies for the optimal proposal versus the bounded approximation are very similar. The geometric mean approximation also provides a reasonable first approximation to the optimal proposal for the examples shown in Figures \ref{fig:gaussian_qplots}--\ref{fig:chi2_qplots}, with improvements in sampling efficiencies compared to the KDE or posterior-approximation proposals (see Table \ref{tab:sampling_efficiencies}).

These three one-dimensional examples have the virtue of being easy to visualise and to show the features of the optimal proposal for Gaussian, skewed and multi-modal posteriors. In the following section we demonstrate that a lower bound on the expected relative improvement of the optimal proposal  improves exponentially on the performance of other kernels.

\subsubsection{Expected improvement as a function of data informativeness and parameter dimensionality}
We can obtain a lower bound on the improvement in sampling efficiency enabled by the optimal proposal by applying the geometric mean approximation to a simple toy problem. In this model, we take the prior and posterior to be $n_\theta$-dimensional multivariate Gaussians, with means in each dimension of $\mu_\pi$ and $\mu_p=0$, and diagonal covariances with variances $\sigma^2_\pi$ and $\sigma^2_p=1$, respectively. We then exploit the fact that  $\omega[q_0]$ can be obtained analytically in this setting to rapidly evaluate the improvement in sampling efficiency over other choices of the proposal as a function of the dimensionality of the parameter space and the informativeness of the data (as measured by both the reduction in volume and the shift in the mean in going from prior to posterior). As the optimal proposal $q^\ast$ outperforms the geometric mean approximation $q_0$, the results in this section can be considered lower bounds on the improvement in performance realized by using the optimal proposal.

\begin{figure}
\centering
 \includegraphics[width=\textwidth]{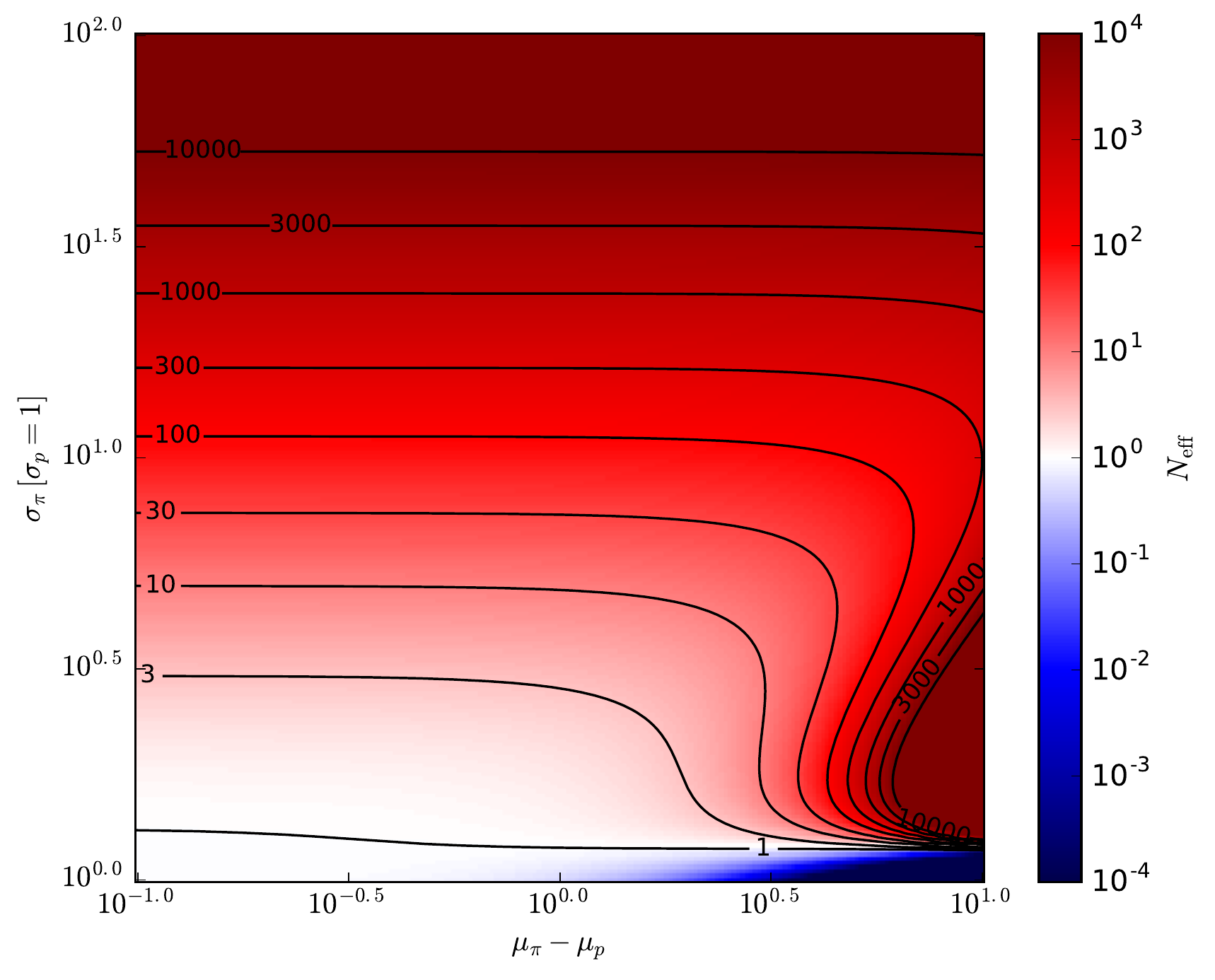}
\caption{Improvement $a$ in the sampling efficiency when sampling from the geometric mean approximation to the optimal proposal rather than the posterior for the case where both prior and posterior are Gaussian. This plot shows the case where $n_\theta=3$ (3-parameters). Even for this modest number of parameters the improvement is large when the data are  informative. While the exact result guarantees $a>1$, the geometric mean approximation gives $a<1$ \textit{(shown in blue)} for the atypical case when the posterior and prior have nearly equal width but are located far from each other.  }\label{fig:approximate_improvement_posterior}
\end{figure}

\begin{figure}
\centering
 \includegraphics[width=\textwidth]{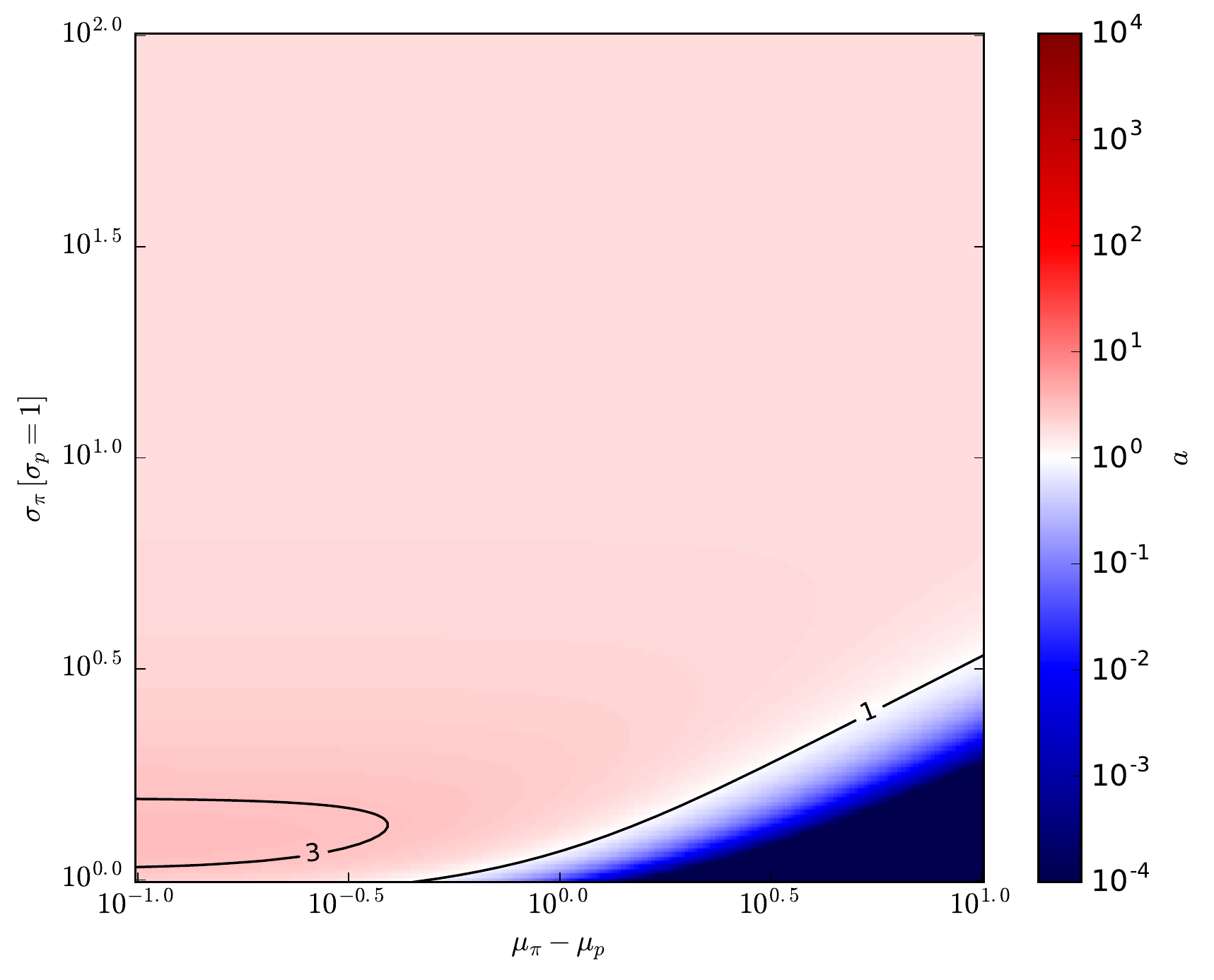}
\caption{Improvement $a$ in the sampling efficiency when drawing proposals from the geometric mean approximation to the optimal proposal rather than the KDE scheme recommended in \cite{Beaumont2009} (with bandwidth equal to twice the posterior variance) for the case where both prior and posterior are Gaussian. This plot shows the case where $n_\theta=3$. While the exact result guarantees $a>1$, the approximation gives $a<1$ \textit{(shown in blue)} for the atypical case when the posterior and prior have nearly equal width but are located far from each other.  }\label{fig:approximate_improvement_literature_3d}
\end{figure}

\begin{figure}
\centering
 \includegraphics[width=\textwidth]{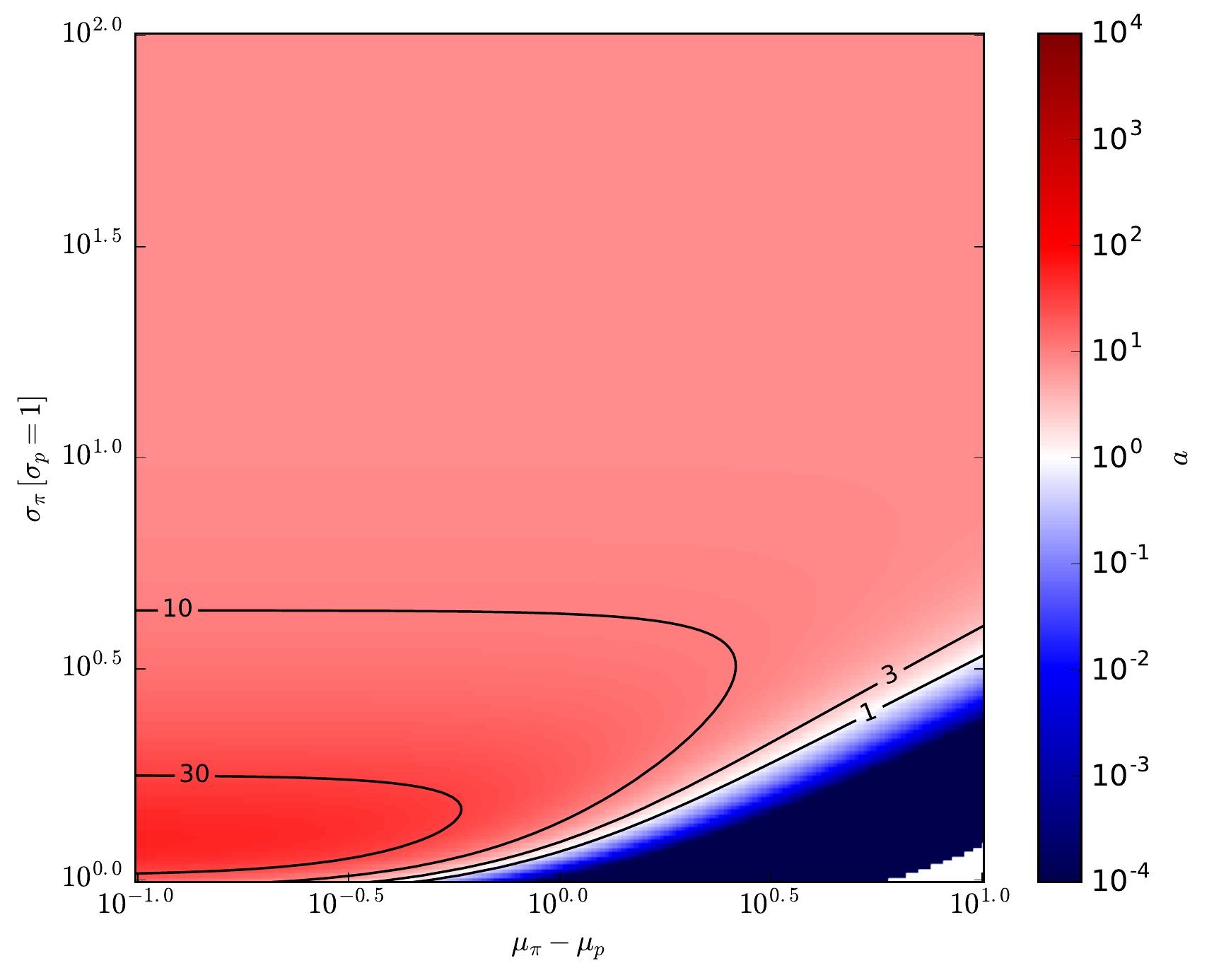}
\caption{Improvement $a$ in the sampling efficiency when drawing proposals from the geometric mean approximation to the optimal proposal rather than the KDE scheme recommended in \cite{Beaumont2009} (with bandwidth equal to twice the posterior variance) for the case where both prior and posterior are Gaussian. This plot shows the case where $n_\theta=10$. While the exact result guarantees $a>1$, the approximation gives $a<1$ \textit{(shown in blue)} for the atypical case when the posterior and prior have nearly equal width but are located far from each other.  }\label{fig:approximate_improvement_literature_10d}
\end{figure}

Figure \ref{fig:approximate_improvement_posterior} compares the geometric mean approximation proposal to using an estimator for the posterior density as the proposal, in the three-dimensional setting. It is clear that proposing samples from the posterior is highly suboptimal, whilst using even the geometric mean approximation to the optimal proposal instead gives improvements that are large, especially when the data are highly informative or surprising compared to prior expectations. Figure \ref{fig:approximate_improvement_literature_3d} compares the geometric mean approximation proposal to the Gaussian KDE proposal scheme recommended in \cite{Beaumont2009} (with bandwidth set to double the estimated posterior variance), again for the case where the posterior and prior are three-dimensional Gaussians. Again, gains in sampling efficiency are expected using the geometric mean approximation as opposed to the KDE scheme. Figure \ref{fig:approximate_improvement_literature_10d} shows the same case but for a 10-parameter rather than three-parameter set-up, showing that the relative improvement from using the optimal kernel quickly becomes larger in higher dimensions.

Inspecting the analytical result for the Gaussian case confirms that the relative improvement in sampling efficiency of the geometric mean approximation scales exponentially with the number of parameters when compared to the other cases we study here. This also implies exponentially scaling improvement with number of parameters for the optimal proposal.


\section{Conclusions}
We have derived an optimal proposal scheme for SMC-ABC algorithms by maximizing the sampling efficiency, defined as the effective number of samples per parameter proposal. This represents the optimal trade-off between obtaining a high acceptance rate, whilst reducing the variance of the importance-weights of the accepted samples, hence increasing the effective sample size for the accepted samples. We derived an implicit form for the optimal proposal that can be solved iteratively, and also two convenient and simple approximations that can be evaluated quickly, provided a density estimator for the posterior is available (based on the accepted samples from the previous SMC population). 

We have shown that the optimal proposal scheme gives a substantial boost in the expected sampling efficiency for a range of statistical models, and the expected gain in sampling efficiency increases with the dimensionality of the problem. The derived results hence provide a guide for choosing optimal proposal densities for SMC-ABC applications, for the small cost of constructing a posterior density estimator at each SMC population iteration.

\bibliographystyle{ba}
\bibliography{optimalABC}

\begin{acknowledgement}
We wish to thank Ethan Anderes  and Tom Charnock for  discussions. This work is supported by the Simons Foundation. Benjamin Wandelt acknowledges support from the Labex Institut Lagrange de Paris (ILP) (reference ANR-10- LABX-63) part of the Idex SUPER, and received financial state aid managed by the Agence Nationale de la Recherche, as part of the programme Investissements d’avenir under the reference ANR-11-IDEX-0004-02.
\end{acknowledgement}
\end{document}